# The convection-diffusion equation for a finite domain with time varying boundaries[1,2,3]

W. J. Golz
Department of Civil and Environmental Engineering, Louisiana State University
Baton Rouge, Louisiana 70803, USA

J. R. Dorroh
Department of Mathematics, Louisiana State University
Baton Rouge, Louisiana 70803, USA

**Abstract**—A solution is developed for a convection-diffusion equation describing chemical transport with sorption, decay, and production. The problem is formulated in a finite domain where the appropriate conservation law yields Robin conditions at the ends. When the input concentration is arbitrary, the problem is underdetermined because of an unknown exit concentration. We resolve this by defining the exit concentration as a solution to a similar diffusion equation which satisfies a Dirichlet condition at the left end of the half line. This problem does not appear to have been solved in the literature, and the resulting representation should be useful for problems of practical interest.

Authors of previous works on problems of this type have eliminated the unknown exit concentration by assuming a continuous concentration at the outflow boundary. This yields a well-posed problem by forcing a homogeneous Neumann exit, widely known as the Danckwerts [1] condition. We provide a solution to the Neumann problem and use it to produce an estimate which demonstrates that the Danckwerts condition implies a zero concentration at the outflow boundary, even for a long flow domain and a large time.

## 1. INTRODUCTION

CLASSICAL SOLUTIONS of the diffusion equation have been catalogued for many of the important problems of heat transfer [2]. The diffusion equation has also been widely employed as a model for chemical reaction processes, and this usually entails the inclusion of lower-order terms that describe convection and reaction. The complications introduced by those additional terms often call for inventive techniques which yield novel and useful representations (e.g., see [3]).

Descriptions of chemical processes and contaminant transport have motivated a large volume of work on the one-dimensional convection-diffusion equation (CDE). Numerical methods have provided solutions to problems satisfying a fairly wide range of conditions (e.g. see [3-5]). However, analytical solutions continue to be highly valued for their inherent simplicity, their capacity to convey qualitative information about the physical problem, and as a verification for numerical models. Many of the existing analytical solutions have been compiled in two large compendiums [6,7].

---

The literature on the mathematics of chemical transport demonstrates that problems posed with a Neumann condition fail to satisfy the relevant conservation equation for concentrations in the interior of the domain [8-10]. The inconsistency with the physical problem arises because the conservation of mass requires Robin conditions at the ends, which leads to a system that is underdetermined due to an unknown exit concentration [1,10]. Danckwerts eliminated the unknown exit concentration by assuming a continuous concentration at the outflow boundary, thus forcing a homogeneous Neumann condition [1,10]. The Danckwerts assumption yields a well-posed problem, but it introduces error into predicted concentrations [6,10]. The magnitude of the error is inversely proportional to time and Peclet number (i.e., the product of fluid velocity and flow length when divided by the diffusion constant, see [6,10]). Hydrogeological problems characteristically exhibit a low Peclet number, and observations at early times are often important. Thus, the error can be significant, especially near an outflow boundary [6,10].

This paper was motivated by a class of physically important problems from contaminant hydrogeology. These problems can be characterized by a diffusion equation with additional terms for drift, decay, and production, and the input may be any bounded continuously differentiable function. The resulting system is composed of an inhomogeneous CDE with time varying Robin conditions. Mass is conserved when the exit concentration is defined as the solution to a CDE with a Dirichlet condition at the left end of the half line. This approach yields the proper formulation in finite geometry thus admitting a representation in an eigenfunction expansion with closed-form eigenvalues.

## 2. FORMULATION AND TRANSFORMATION

The problems that will be considered in this paper take the form

$$RC_t = DC_{xx} - vC_x - \mu C + \gamma, \qquad 0 < x < \ell, \quad t \in \mathbb{R} \tag{2.1}$$

where $C = C(x, t)$ denotes a concentration, with x the longitudinal distance and t the time. The constants R, D, v, $\mu$, and $\gamma$ describe linear sorption, diffusion, longitudinal fluid velocity, decay, and production, respectively. Equation (2.1) will satisfy the auxiliary conditions

$$C(x, t_0) = \phi(x), \qquad 0 < x < \ell \tag{2.2}$$

$$vC(0, t) - DC_x(0, t) = vg(t), \qquad t > t_0 \tag{2.3}$$

$$vC(\ell, t) - DC_x(\ell, t) = C_E, \qquad t > t_0 \tag{2.4}$$

where $\phi(x), g(t)$ are bounded continuously differentiable functions, and $C_E$ is the exit concentration which, for now, we regard as an experimentally measured quantity.

To transform the CDE into standard form for an inhomogeneous diffusion equation, and to provide homogeneous boundary conditions, we introduce the change of dependent variable

$$C(x, t; r, s) = \left(w(x, t) + e^{st} H(x, t)\right) e^{rx - st} \tag{2.5}$$

If we choose the parameters as

$$r = \frac{v}{2D}, \qquad s = \frac{1}{R}\left(\frac{v^2}{4D} + \mu\right), \qquad D, R \neq 0 \tag{2.6}$$

then (2.1) may be written as the inhomogeneous diffusion equation

$$w_t = \frac{D}{R} w_{xx} + e^{st} F(x, t), \qquad F(x, t) = \frac{\gamma}{R} e^{-rx} - H_t + \frac{D}{R} H_{xx} \tag{2.7}$$

and the initial condition (2.2) will become

$$w(x, t_0) = e^{st_0}\left(e^{-rx}\phi(x) - H(x, t_0)\right) \tag{2.8}$$

If we then define

$$H(x, t) = \left(1 + \cos\frac{\pi x}{\ell}\right)g(t) + \left(1 - \cos\frac{\pi x}{\ell}\right)e^{-\ell x}C_E \tag{2.9}$$

the boundary conditions (2.3), (2.4) will be homogeneous:

$$w_x(0, t) - r w(0, t) = 0 \tag{2.10}$$

$$w_x(\ell, t) - r w(\ell, t) = 0 \tag{2.11}$$

## 3. EIGENFUNCTION EXPANSION

To obtain our representation, we only require the separated solution $w(x, t) = \varphi(x)T(t)$ to be justified when $F(x, t) = 0$. This leads to a regular Sturm-Liouville problem with the positive eigenvalues and related eigenfunctions

$$\lambda_n = \frac{n^2 \pi^2}{\ell^2}, \qquad \varphi_n(x) = \cos\left(x\sqrt{\lambda_n}\right) + \frac{r}{\sqrt{\lambda_n}}\sin\left(x\sqrt{\lambda_n}\right), \qquad n = 1, 2, 3, \cdots \tag{3.1}$$

where each $\lambda_n$ defines a single-dimensional eigenspace. There is also one negative eigenvalue with a single associated eigenfunction:

$$\lambda_n = -r^2, \qquad \varphi_n(x) = e^{rx}, \qquad n = 0 \tag{3.2}$$

The negative eigenvalue arises from the transformations, so its related eigenfunction does not entail the usual meaning of a source term which would allow the solution to become arbitrarily large as t increased. In fact, we will exhibit a large-t solution that is bounded.

The normalized eigenfunctions of a Sturm-Liuoville system clearly compose an orthonormal basis (e.g., see [11], sec. 7.5), from which it follows that (2.7) will have a solution of the form

$$w(x, t) = \sum_{n=M}^{\infty} \varphi_n(x) T_n(t), \qquad M = 0 \tag{3.3}$$

where

$$T_n(t) = \frac{e^{-\frac{D}{R}\lambda_n t}}{\int_0^\ell \varphi_n^2(x)\,dx}\left(\int_{t_0}^t e^{\left(s + \frac{D}{R}\lambda_n\right)\tau}\int_0^\ell F(x, \tau)\varphi_n(x)\,dx\,d\tau \right.$$

$$\left. + e^{\left(s + \frac{D}{R}\lambda_n\right)t_0}\int_0^\ell \left(e^{-rx}\phi(x) - H(x, t_0)\right)\varphi_n(x)\,dx\right) \tag{3.4}$$

To obtain (3.4), the orthogonality of the eigenfunctions and the symmetry of the boundary conditions on $w$ and $\varphi$ are required to develop an ordinary differential equation in T that satisfies the initial condition in w (for a discussion of this technique, see [12], chap. 9).

We now wish to demonstrate that the representation of $C(x,t)$ will be bounded for all t. To establish this, we fix t in (3.4) and then note that the second term is bounded and will therefore vanish as $t_0 \to -\infty$. Then referring to (2.5), a large-t solution to (2.1)-(2.4) will be

$$C_{bv}(x,t) = \sum_{n=M}^{\infty} \frac{\varphi_n(x) e^{rx - \left(s + \frac{D}{R}\lambda_n\right)t}}{\int_0^{\ell} \varphi_n^2(x)\,dx} \int_{-\infty}^{t} e^{\left(s + \frac{D}{R}\lambda_n\right)\tau} \int_0^{\ell} F(x,\tau)\varphi_n(x)\,dx\,d\tau + e^{rx} H(x,t) \quad (3.5)$$

which is obviously nonincreasing for t. Note that (3.5) is equivalent to a solution of the pure boundary value problem (2.1), (2.3), (2.4) where it is valid for all $t \in \mathbb{R}$ (for a similar result involving a problem on the half line, see [3]).

## 4. EXIT CONCENTRATION

In the event that the exit concentration is not empirically observed, we may determine it. First, set

$$vC_F(x,t) = vC(x,t) - DC_x(x,t), \qquad 0 \leq x < \infty, \quad t \in \mathbb{R} \quad (4.1)$$

It can be verified by direct substitution that $C_F$ satisfies a differential equation identical to (2.1) (e.g., see [8,9]). It then follows that $C_E = C_F(\ell, t)$, where $C_F(x, t)$ is a bounded solution to (2.1) with the auxiliary conditions

$$vC_F(x,t_0) = v\phi(x) - D\phi_x(x), \qquad 0 < x < \infty \quad (4.2)$$

$$C_F(0,t) = g(t), \qquad t > t_0 \quad (4.3)$$

The change of variable

$$C_F(x,t) = u(x,t) e^{rx - st} + \frac{\gamma}{\mu} \quad (4.4)$$

will, with reference to (2.1), yield the canonical form of the diffusion equation

$$u_t = \frac{D}{R} u_{xx} \quad (4.5)$$

and (4.2), (4.3) will become

$$u(x, t_0) = \Phi(x) \quad (4.6)$$

$$u(0, t) = G(t) \quad (4.7)$$

where

$$\Phi(x) = \left(\phi(x) - \frac{D}{v}\phi_x(x) - \frac{\gamma}{\mu}\right) e^{-rx + st_0}, \qquad G(t) = \left(g(t) - \frac{\gamma}{\mu}\right) e^{st} \quad (4.8)$$

A solution to (4.5)-(4.7) is provided by the familiar integral representation

$$u(x,t) = \int_0^{\infty} K(x - \zeta, (D/R)(t - t_0)) - K(x + \zeta, (D/R)(t - t_0)) \Phi(\zeta)\,d\zeta \\ - \frac{2D}{R} \int_{t_0}^{t} K_x(x, (D/R)(t - \tau)) G(\tau)\,d\tau \quad (4.9)$$

where the first term may be obtained from an odd extension of $\Phi$ and the last term from the method of Duhamel (e.g., see [2], chap. 4). The kernel and its derivative are then

$$K(x,t) = \frac{1}{\sqrt{4\pi t}} e^{-x^2/4t}, \qquad K_x(x,t) = -\frac{x}{2t} K(x,t) \quad (4.10)$$

and may be obtained in a straightforward manner from a Fourier-transform solution of the fundamental initial-value problem (e.g., see [2], chap. 3) or from a Green's function (e.g., see [3]).

## 5.   AN ESTIMATE OF THE DANCKWERTS ERROR

In this section, we provide a solution to the Danckwerts problem (i.e., the Robin exit condition is replaced with a Neumann condition). The resulting representation is then used to produce an estimate of the Danckwerts error.

First we define a new concentration $C_D = C_D(x,t)$ which will satisfy (2.1)-(2.3), but (2.4) will be replaced by the Neumann condition

$$C_{Dx}(\ell, t) = 0, \qquad t > t_0 \tag{5.1}$$

To see that (5.1) is not generally valid, take $C_{Dt} = \gamma = 0 \ni C_D = C_D(x)$ in (2.1) which yields the problem considered by Danckwerts [1]. Now note that $C_{Dx}(\ell) = 0 \Rightarrow C_D(\ell) = 0$.

The change of variable (2.5) and parameters (2.6) remain the same while the CDE (2.7) and the initial condition (2.8) must respect the new definition

$$H_D(x,t) = \left(1 + \cos\frac{\pi x}{\ell}\right) g(t) \tag{5.2}$$

which replaces (2.9). The entrance condition (2.10) for $w_D = w_D(x,t)$ will remain unchanged, but the exit condition (2.11) will now be replaced by

$$w_{Dx}(\ell, t) + r\, w_D(\ell, t) = 0 \tag{5.3}$$

The separated solution $w_D(x,t) = \varphi_D(x) T_D(t)$ for the homogeneous problem in $w_D$ will now yield a Sturm-Liouville problem with strictly positive eigenvalues. The eigenvalues will no longer be given in the closed form of (3.1) but will instead arise from the sequential intersections $z_1 = z_2$, where

$$z_1 = \tan\left(\ell \sqrt{\lambda_D}\right), \qquad z_2 = \frac{2r\sqrt{\lambda_D}}{\lambda_D - r^2} \tag{5.4}$$

Information about how the $\lambda_D$ approach their asymptotic value is useful for evaluating the error estimate. First, observe that

$$\frac{n^2 \pi^2}{\ell^2} < \lambda_{Dn} < \frac{(n+1)^2 \pi^2}{\ell^2} \tag{5.5}$$

so for each n, $\lambda_n < \lambda_{Dn}$. Then, the $\lambda_{Dn}$ will have the asymptotic value

$$\lim_{\sqrt{\lambda_D} \to \infty} z_2 = \lim_{\sqrt{\lambda_D} \to \infty} \frac{2r\sqrt{\lambda_D}}{\lambda_D - r^2} = \lim_{\sqrt{\lambda_D} \to \infty} \frac{r}{\sqrt{\lambda_D}} = 0 \Rightarrow \lim_{n \to \infty} \lambda_{Dn} = \frac{n^2 \pi^2}{\ell^2} \tag{5.6}$$

The eigenfunctions $\varphi_{Dn}(x)$ will be defined as in (3.1) but with their eigenvalues now given by (5.4). If we let the summation begin at $M = 1$, a solution for $w_D(x,t)$ will be provided by (3.3) where the $T_{Dn}(t)$ come from (3.4), and a large-t solution will be given by $C_{Dbv}(x,t)$ as in (3.5).

Since C and $C_D$ differ most near $x = \ell$, our definition of the error is the natural one

$$E_D(t) = |C(\ell,t) - C_D(\ell,t)| \tag{5.7}$$

If we allow that (5.7) is, by hypothesis, nonincreasing with t, it then follows that

$$E_D(t) \geq |C_{bv}(L,t) - C_{Dbv}(L,t)| \geq |\gamma/\mu| \tag{5.8}$$

The verification of (5.8) is left to the reader, but it is easy to see that $|\lambda/\mu|$ will be obtained from $|C_{bv}(\ell,t) - C_{Dbv}(\ell,t)|$ at a sufficiently large constant length L since $\lambda_n, \lambda_{Dn} \to 0$ for the first few eigenvalues and similar terms in the large-t solutions $C_{bv}(\ell,t), C_{Dbv}(\ell,t)$ approach one another quite rapidly as n increases.

We now refer to (2.1), where it is clear that $|\lambda/\mu|$ is the amount that C(x,t) will differ from zero at a very large distance. Thus, the Neumann condition implies a zero concentration at the outflow boundary, just as in the original problem considered by Danckwerts [1].